# A rate-induced tipping in the Pearson diffusion

Hidekazu Yoshioka[1][0000-0002-5293-3246]

[1] Japan Advanced Institute of Science and Technology, 1-1 Asahidai, Nomi, 923-1292, Japan
yoshih@jaist.ac.jp

**Abstract**. Rate-induced tipping is an instability that occurs in a system when its time-dependent rate parameter becomes larger than a threshold value. We investigate a Pearson diffusion process, a diffusion process having solutions staying in a bounded domain under certain conditions, whose noise-free limit experiences a rate-induced tipping such that solutions escape from the domain in a finite time. We show that the existence of noise leads to faster escapement from the domain.

## 1 Introduction

### 1.1 Pearson diffusion

Analyzing stability of a dynamical system is a basic strategy towards understanding its properties. Stability of a dynamical system at boundary is crucial especially when solutions to the system are confined in a prescribed domain; a situation where the escape of solutions from the domain represents some instability phenomenon. For stochastic dynamical systems, such an analysis often reduces to studying boundary behavior of a solution to a stochastic differential equation (SDE) [1-3].

The Pearson diffusion, also called Jacobi diffusion, presented below is the simplest SDE whose solution is almost surely (a.s.) confined in the closed unit interval $D = [0,1]$ under certain conditions (Chapter 6 in Alfonsi [4]):

$$\mathrm{d}X_t = \underbrace{(r - X_t)\mathrm{d}t}_{\text{Drift}} + \underbrace{\sigma\sqrt{X_t(1-X_t)}\mathrm{d}B_t}_{\text{Diffusion}}, \ 0 < t < \tau \qquad (1)$$

subject to an initial condition $X_0 \in D_{\text{int}} = (0,1)$. Here, $t \geq 0$ is time, $r > 0$ is source rate, $\sigma > 0$ is volatility, and $B = (B_t)_{t \geq 0}$ is a 1-D standard Brownian motion (we will work with a natural filtration generated by $B$), and $\tau$ is the stopping time defined by

$$\tau = \inf\{t > 0 \mid X_t \geq 1 \text{ or } X_t \leq 0\}. \qquad (2)$$

Namely, $\tau$ is the escape time of the process $X$ from $D$. The diffusion term in the SDE (1) is understood in the Itô's sense.

In (1), the drift term represents mean reversion of the process $X$ to the level $r$, and the diffusion term represents stochastic fluctuations that randomly drives $X$. Due to the forms of the drift and diffusion terms, the boundary behavior and hence stability of solutions to the SDE (1) is found as follows (Theorem 6.1.1 and Proposition 6.1.2 in Alfonsi [4]):

- ✓ If $0 < r \leq 1$, then the SDE (1) admits a unique path-wise solution that is continuous in time. Moreover, the solution is a.s. confined in $D$.
- ✓ The solution to the SDE (1) is a.s. confined in $D_{\text{int}}$ ($\tau = +\infty$) if and only if $2r \geq \sigma^2$ and $2(1-r) \geq \sigma^2$.

The condition $0 < r \leq 1$ in the first bullet point is necessary for $2(1-r) \geq \sigma^2$ in the second one. More specifically, if $2(1-r) < \sigma^2$, then the solution to the SDE (1) has an a.s. finite $\tau$ with $X_\tau = 1$.

### 1.2  Our model

Boundary stability of the Pearson diffusion has thus been well-studied. Particularly, the size of the parameter $r$ is critical for the stability. By contrast, the following time-dependent version where the source $r$ is gradually decreasing has not been studied in the literature to the best of the author's knowledge:

$$dX_t = (Y_t - X_t)dt + \sigma\sqrt{X_t(1-X_t)}dB_t, \ 0 < t < \tau \tag{3}$$

and ($\tau$ is again given by (2))

$$X_t = 1 \text{ if } t \geq \tau \tag{4}$$

coupled with

$$dY_t = f(Y_t)dt, \ 0 < t < \tau \tag{5}$$

subject to an initial condition $(X_0, Y_0) \in D_{\text{int}} \times \mathbb{R}$, where $Y = (Y_t)_{t \geq 0}$ is a deterministic process that represents time-dependent source rate, and $f : \mathbb{R} \to \mathbb{R}$ is a function such that the ordinary differential equation (ODE in short) (5) admits a unique smooth solution.

The objective of this study is to quantify stability of the novel system (3)-(5). A dynamical system (the SDE (3) and the relationship (4)) driven by a temporally decreasing or increasing external force (the ODE (5)) has been studied in the context of rate-induced tipping [5]. Here, rate-induced tipping is an instability that occurs in a system when its rate parameter becomes larger than a threshold value. In our case, the tipping and hence the escape of the process $X$ from the domain $D$ possibly occurs if $2Y_t < \sigma^2$ or $2(1-Y_t) < \sigma^2$. Particularly, the escape from the upper boundary of $D$ possibly occurs if $2(1-Y_t) < \sigma^2$. Because the system (3)-(5) is stochastic, the stabil-ity should be evaluated using some probability measure that characterizes the tipping phenomenon.

Before starting the analysis of the system (3)-(5), we discuss when such a system may arise. Yoshioka [6] modelled the traveler's demand (i.e., normalized arrival intensity of travelers at a tourism spot) by using a Pearson diffusion process. In this case, the source rate is the attractiveness of the tourism spot, and a too large value of $X$ is assumed to trigger an overtourism. Then, a time-dependent source rate would represent some tourism restriction by a social planner so that the tourism becomes sustainable (i.e., $X$ is a.s. confined in $D_{\text{int}}$ with a high probability).

## 2   Problem setting

### 2.1   The model

The system studied in this paper is given by (3)-(5) where the coefficient $f$ is set as the convex quadratic function

$$f(y) = Ry\left(1 - \frac{y}{1-\Delta}\right), \quad y \in \mathbb{R} \tag{6}$$

with parameters $R > 0$ and $\Delta \in (0,1)$ subject to a prescribed initial condition $(X_0, Y_0) \in D_{\text{int}} \times \{1+\Delta\}$. A rationale behind the choice of the coefficient (6) along with the initial condition $Y_0 = 1+\Delta$ is that the source rate initially violates the condition of boundary stability $2(1-Y_t) \geq \sigma^2$, while this condition is satisfied after some later time because the solution to the ODE (5) gradually decreases from $1+\Delta$ to $Y_0 = 1-\Delta$ as the time elapses. We therefore infer that the escape of the process $X$ from the domain $D$ through the upper boundary possibly occur when $2(1-Y_t) < \sigma^2$, namely during a time interval close to the initial time 0. The escape probability $p$ of our system is defined as

$$p = \mathbb{P}(X_t = 1 \text{ at some } t > 0). \tag{7}$$

The ODE (5) with the coefficient (6) is inspired from the literature of rate-induced tipping [5]. Indeed, the ODE (5) has a solution having a sigmoidal shape that is controlled by the solo parameter $R$: a larger $R$ represents faster variation of $Y$ from $1+\Delta$ to $1-\Delta$. In this view, we infer that a larger value of $R$ leads to a smaller probability of the escapement of $X$ from the upper boundary of $D$.

Finally, for the deterministic case ($\sigma = 0$), it has been numerically found that there exists some $R = R_c$ such that $X_t$ is always strictly smaller than 1 if $R > R_c$ and there exists some $t = t_c$ such that $X_{t_c} = 1$ if $R \leq R_c$. In this view, the deterministic system experiences a rate-induced tipping at $R = R_c$.

## 2.2 Monte Carlo discretization

We numerically discretize the system (3)-(5) along with (6) because we do not have closed-form solutions to this system. We use a Monte Carlo method based on the common Euler–Maruyama method that explicitly discretize the system in time. We also equip the discretized system with a truncation so that numerical solutions inherit the property "$X_t = 1$ if $t \geq \tau$" in (3): for $n = 0, 1, 2, ...$

$$X_{(n+1)h} = \begin{pmatrix} X_{nh} + (Y_{nh} - X_{nh})h + \sigma\sqrt{X_{nh}(1-X_{nh})} \times \sqrt{h}W_n & (X_{nh} < 1) \\ X_{nh} & (X_{nh} \geq 1) \end{pmatrix}, \quad (8)$$

$$Y_{(n+1)h} = Y_{nh} + RY_{nh}\left(1 - \frac{Y_{nh}}{1-\Delta}\right)h \quad (9)$$

subject to an initial condition $(X_0, Y_0) \in D_{int} \times \{1+\Delta\}$. Here, $h > 0$ is a sufficiently small time increment, the subscript $nh$ represents the quantity approximated at time $t = nh$, and $W_n$ is a standard normal random variable with mean 0 and variance 1 where each $W_n$ being assumed to be mutually independent. For each computational case presented in this paper, we use the increment $h = 10/200,000$ and compute the escape probability by using 200,000 sample paths. These computational resolutions have been found to be sufficiently fine by preliminary computational experiments.

## 3 Computational investigation

### 3.1 Sample path, average, standard deviation

We fix the parameter values $X_0 = 0.1$ and $\Delta = 0.5$. **Fig. 1** shows the computed sample paths of $X$ where different colors represent different sample paths, demonstrating that a larger value of $R$ yields the process $X$ confined in $D_{int}$. **Fig. 2** shows the computed average (Ave) and standard deviation (Std) of the corresponding computational cases. These statistics reflect the behavior of sample paths plotted in **Fig. 1**.

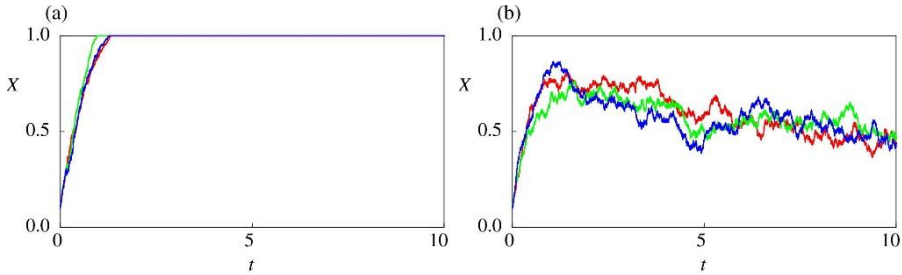

**Fig. 1.** Computed sample paths of $X$ where different colors represent different sample paths: (a) $(\sigma, R) = (0.2, 0.1)$ and (b) $(\sigma, R) = (0.2, 0.5)$.

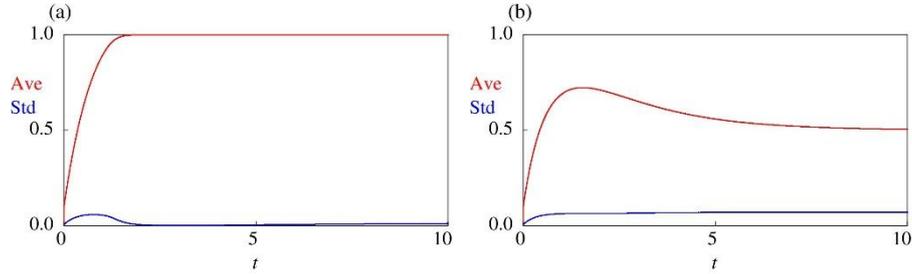

**Fig. 2.** Computed Ave and Std of: (a) $(\sigma, R) = (0.2, 0.1)$ and (b) $(\sigma, R) = (0.2, 0.5)$.

### 3.2 Escape probability and hitting time

**Table 1** summarizes the computed escape probabilities $p$ for different couples of $(\sigma, R)$, showing that $p$ is decreasing with respect to $\sigma$ and $R$. The stability of the system in terms of the escape probability $p$ is decreasing with respect to $R$ as expected. Increasing the volatility $\sigma$, which is the noise amplitude, increases the escape probability $p$. The decrease of $p$ with respect to $R$ is more gradual for larger $\sigma$.

To deeper understand the escapement phenomenon, **Fig. 3** shows histograms of the computed first hitting times $\tau$ of the process $X$ at the upper boundary of the domain $D$. Notice that now $\tau = \inf\{t > 0 | X_t = 1\}$. According to **Fig. 3**, the histograms are unimodal and are flatten as $\sigma$ increases.

**Table 1.** Escape probabilities $p$ for different couples of $(\sigma, R)$ represented in %.

| $R \backslash \sigma$ | 0.1 | 0.2 | 0.4 | 0.8 |
|---|---|---|---|---|
| 0.1 | 100.00 | 99.89 | 98.83 | 99.70 |
| 0.2 | 0.00 | 0.58 | 28.01 | 87.68 |
| 0.3 | 0.00 | 0.00 | 1.12 | 61.96 |
| 0.4 | 0.00 | 0.00 | 0.02 | 39.16 |
| 0.5 | 0.00 | 0.00 | 0.00 | 23.78 |

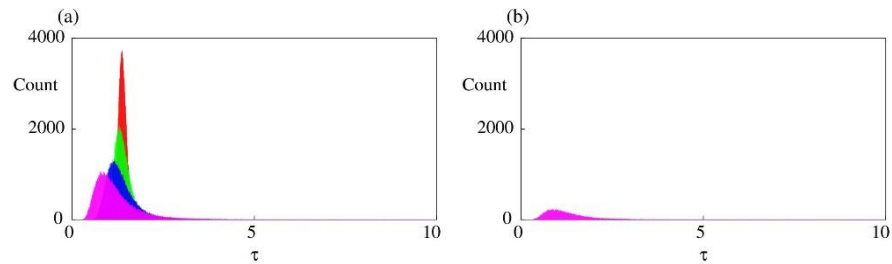

**Fig. 3.** Histograms of the computed first hitting times $\tau$ of the process $X$ at the upper boundary of the domain $D$: (a) $R = 0.1$ and (b) $R = 0.5$, and $\sigma = 0.1$ (red), $\sigma = 0.2$ (green), $\sigma = 0.4$ (blue), $\sigma = 0.8$ (magenta).

## 4   Conclusion

We studied a stochastic system based on a Pearson diffusion driven by a time-dependent decreasing source rate. The value and decreasing speed of the source rate was found to affect the boundary stability of the solution to the system. The escape probability as a measure of sustainability of the system was computationally quantified.

The system studied in this paper does not have any mean filed effects such that behavior of a sample path of the system is affected by the other paths. This kind of effect crucial when modeling social phenomena such as tourist dynamics. We will study this issue from the standpoint of some mean filed game where the system will be described by SDEs having law-dependent coefficients.

**Acknowledgements**
This study was supported by the Japanese Society for the Promotion of Science (KAKENHI, No. 22K14441), Nippon Life Insurance Foundation (Environmental problems research grant for young researchers, No. 24), and Japan Science and Technology Agency (PRESTRO, No. 24021762).